\input amstex
\documentstyle{amsppt}
\magnification=1200
\hsize = 16 true cm
\vsize = 22 true cm

\document
\overfullrule=0pt \baselineskip=1.2 \normalbaselineskip
\NoRunningHeads

\topmatter
\title
 A Hybrid Volterra-type Equation with Two Types of Impulses
\endtitle

\author
${}^*$S.A. Belbas and ${}^{**}$Jong Seo Park
\bigskip
${}^{*}$Mathematics Department, University of Alabama, Box 870350,
Tuscaloosa, AL 35487-0350, USA
${}^{**}$Department of Mathematics Education, Chinju \\
National University of Education, Jinju,
660-756, Korea\\
\endauthor
\thanks
Keywords and Phrases :  Volterra integral equation, Impulsive
equations, Fixed point theorems.
\endthanks
\thanks
AMS 2000 Mathematics Subject Classification : 45D05, 45L05, 45G10.
\endthanks

\NoRunningHeads
\endtopmatter

\document
\bigskip
\heading Abstract
\endheading
We formulate and analyze a hybrid system model that involves
Volterra integral operators with multiple integrals and two types
of impulsive terms. We give a constructive proof, via an iteration
method, of existence and uniqueness of solutions.

\bigskip
\heading I. Introduction \endheading
\bigskip
The class of Volterra integral equations has been traditionally
used to model the behavior of systems with memory. The question of
existence and uniqueness of solutions of Volterra integral
equations is carried out by applying appropriate fixed point
theorems; in addition to the classical works $[4, 8]$, we mention
the papers $[2, 9, 10]$ that contain techniques related to the
material of the present paper.

In this paper, we formulate and analyze a novel class of equations
that contain some of the characteristics of Volterra  integral
equations, but also fall into the general category of hybrid
systems theory.

The discipline of hybrid systems aims to formulate and analyze
models of systems that combine continuous and discrete features,
as well as to solve associated problems of control theory and game
theory. In the area of ordinary differential equations, the
standard model of hybrid systems is the class of impulsive
differential equations [1]. The topic of Volterra integral
equations with impulses has received considerable attention in the
recent research literature$[5, 6, 7, 11]$. To the best of our
knowledge, the paper$[3]$ is the first to contain a general
Volterra type term describing the impulses; previous papers
contained particular types of impulses.

In this paper, we have set out to explore the inclusion of
multiple-integral Volterra terms and multiple types of discrete
impulsive terms, with two qualitatively different types of
impulses; the various ingredients can be combined in several
different ways. It is expected that the present work will form the
basis for further developments in the area of systems described by
Volterra equations with multiple integral terms and also multiple
integral-sum terms. In order to make these concepts clear at this
stage, we shall describe the various types of integral equations
that can arise as extensions and generalizations of the classical
Volterra model.

The standard (nonlinear) Volterra equation of the second kind for
a scalar-valued unknown function $y(t)$ is
$$
y(t)=y_0 (t)+\int _0^t f(t,s,y(s))ds  \leqno (1.1)
$$
It is well known that general smooth phenomena with memory are not
always modelled by single integral operators of the type $\int
_0^t f(t,s,y(s))ds$, but rather they involve series of multiple
integrals. The multiple integral series analogue of $(1.1)$ is
$$
\aligned y(t)&=y_0 (t) +\sum _{n=1}^\infty \frac{1}{n!} \int _0^t
\cdots \int _0^t \int _0^t f_n (t,s_1 , s_2 , \cdots \cdots , s_n
, \\
&\qquad \qquad \qquad \quad \quad y(s_1 ),y(s_2 ), \cdots , y(s_n
))ds_1 ds_2 \cdots ds_n
\endaligned \leqno (1.2)
$$
and a continuous solution is sought over $0 \leq t \leq T$ .

If the functions $f_n$ are continuous in all their arguments,
satisfy bounds
$$
\vert f(t, s_1 ,\cdots , s_n ,x_1 , x_2 ,\cdots
,x_n )\vert \leq C_n
$$
and Lipschitz conditions
$$
\vert f_n (t, s_1 ,\cdots , s_n ,x_1 , x_2 ,\cdots ,x_n )-f_n (t,
s_1 ,\cdots , s_n ,\xi_1 , \xi_2 ,\cdots ,\xi_n )\vert \leq L_n
\sum _{i=1}^n \vert x_i -\xi_i \vert
$$
with the series $\sum_{i=1}^n M_n \frac{t^n}{n!}$,
$\sum_{n=0}^\infty L_{n+1} \frac{t^n}{n!}$ having radius of
convergence $R>T$, then the existence and uniqueness of a solution
of $(1.2)$ can be shown, in a constructive manner, by showing that
the operator S, defined by
$$
\aligned (Sx)(t):&=y_0 (t) +\sum _{n=1}^\infty \frac{1}{n!} \int
_0^t \cdots \int _0^t \int _0^t f_n (t,s_1 , s_2 , \cdots \cdots ,
s_n , \\
&\qquad \qquad \qquad \quad \quad x(s_1 ),x(s_2 ), \cdots , x(s_n
))ds_1 ds_2 \cdots ds_n
\endaligned \leqno (1.3)
$$
is a contraction on $C(0, T; R)$(the space of continuous
real-valued functions on $[0, T]$) in the norm
$$
\Vert x\Vert:=\max_{0\leq t\leq T} e^{-\mu t}\vert x(t)\vert .
\leqno (1.4)
$$
If $\mu$ is sufficiently large. The proof of the contraction
property relies on the estimate
$$
\aligned
&e^{-\mu t} \vert x(t)-\xi (t)\vert\\
&\leq e^{-\mu t} \sum_{n=1}^\infty \frac{1}{n!} \int _0^t \cdots
  \int _0^t \int _0^t \vert f_n (t,s_1 , s_2 , \cdots ,s_n , x(s_1
  ),x(s_2 ), \cdots , x(s_n ))\\
&\quad \quad \quad \quad \quad \quad \quad -f_n (t,s_1 , s_2 ,
\cdots
   ,s_n ,
  \xi(s_1 ),\xi(s_2 ), \cdots , \xi(s_n ))\vert ds_1 ds_2 \cdots
   ds_n\\
&\leq e^{-\mu t} \sum_{n=1}^\infty \frac{1}{n!} \int _0^t \cdots
  \int _0^t \int _0^t \sum_{j=1}^n L_n \vert x(s_j )-\xi (s_j
  )\vert\\
&\leq e^{-\mu t} \sum_{n=1}^\infty \frac{1}{n!} \int _0^t \cdots
  \int _0^t \int _0^t \sum_{j=1}^n e^{\mu s_j}L_n \Vert x-\xi \Vert _\mu \\
&=\frac{1-e^{-\mu t}}{\mu} \sum_{n=1}^\infty
  \frac{t^{n-1}}{n!}nL_n \Vert x-\xi \Vert _\mu\\
&=\frac{1-e^{-\mu t}}{\mu} \sum_{m=0}^\infty
  \frac{t^{m}}{m!}L_{m+1} \Vert x-\xi \Vert _\mu\\
&\leq \frac{1-e^{-\mu T}}{\mu} \sum_{m=0}^\infty
  \frac{T^{m}}{m!}L_{m+1} \Vert x-\xi \Vert _\mu
\endaligned \leqno (1.5)
$$
and the contraction property follows from the convergence of
$\sum_{m=0}^\infty \frac{T^m}{m!}L_{m+1}$ and the fact that
$\lim_{\mu \longrightarrow \infty}\frac{1-e^{-\mu T}}{\mu} =0$.
\bigskip
It is plain that we may assume, without loss of generality, that
each function $f_n$ is symmetric with respect to permutations of
the symbols $(s_i ,x_i )$, $i=1, 2, \cdots , n$. i.e. $f_n
(t,s_{\sigma(1)} , \cdots ,s_{\sigma (n)} ,x_{\sigma (1)} ,\cdots
, x_{\sigma (n)} )=f_n (t,s_1 , \cdots ,s_n ,x_1 ,\cdots , x_n )$
for every permutation $\sigma$ of $(1, 2,... , n)$, otherwise we
could replace each $f_n$ by its symmetrization
$$
\widetilde{f_n} (t,s_1 , \cdots ,s_n ,x_1 ,\cdots , x_n
):=\frac{1}{n!} \sum_{\sigma \in \prod_n}f_n (t,s_{\sigma(1)} ,
\cdots ,s_{\sigma (n)} ,x_{\sigma (1)} ,\cdots , x_{\sigma (n)} ),
\leqno (1.6)
$$
where $\prod_n$ is the set of all permutations of $(1, 2,... , n)$
, and we would have
$$
\aligned
&\int _0^t \cdots \int _0^t \int _0^t f_n (t,s_1 , s_2 ,
  \cdots ,s_n ,x(s_1 ),x(s_2 ), \cdots , x(s_n ))ds_1 ds_2 \cdots
  ds_n\\
&=\int _0^t \cdots \int _0^t \int _0^t \widetilde{f_n} (t,s_1 ,
s_2 , \cdots ,s_n ,x(s_1 ),x(s_2 ), \cdots , x(s_n ))ds_1 ds_2
\cdots ds_n .
\endaligned \leqno (1.7)
$$
Under the condition of symmetry of each $f_n$, $(1.2)$ can be
written as
$$
\aligned y(t)&=y_0 (t) +\sum _{n=1}^\infty \int _{s_n =0}^t \int
_{s_{n-1} =0}^{s_n} \cdots \int _{s_1 =0}^{s_2} f_n (t,s_1 , s_2 ,
\cdots , s_n
, \\
&\qquad \qquad \qquad \quad \quad y(s_1 ),y(s_2 ), \cdots , y(s_n
))ds_1 \cdots ds_{n-1} ds_n .
\endaligned \leqno (1.8)
$$
Another extension of the standard Volterra equation $(1.1)$ has
been introduced in $[3]$. This extension involves impulsive terms,
and the impulsive Volterra equation becomes
$$
y(t)=y_0 (t)+\int_0^t f(t,s,y(s))ds+\sum_{i:\tau_i <t} g(t,\tau_i
,y(\tau_i )). \leqno (1.9)
$$
(The problems in $[3]$ included controlled equations of the type
$(1.9)$.)

The time-instants $\tau_i$ in $(1.9)$ are the impulsive times. The
unknown function $y$ will have jump discontinuities at these
instants; however, the effect of the jumps in $y$ are more general
than in the case of impulsive ordinary differential equations. A
crucial feature of an impulsive Volterra equation of the type
$(1.9)$ is that a constructive proof of existence and uniqueness
of solutions requires the use of a two-dimensional vector-valued
metric and a concept of contractions with a $2\times 2$ matrix,
instead of a scalar, contraction "coefficient". We refer to $[3]$
for the details.

In this paper, we have further explored the interrelationships
between continuous and discrete components of systems with memory.
To this effect, we have considered Volterra equations with
integral terms up to order two, discrete terms up to order two,
and mixed discrete and continuous terms up to order two, with two
kinds of impulsive effects, impulses of order 1, at impulsive
times denoted by $\tau_i$ , and impulses of order 2, at variable
times denoted by $\sigma_i (t)$. The detailed formulation of this
model, together with explanations of the terminology of
higher-order impulses, is presented in the next section.

The constructive proof of existence and uniqueness of solutions
for our second-order full impulsive model requires the use of a
three-dimensional vector-valued metric, and proof of contraction
property with the role of contraction coefficient played by a
matrix.
\bigskip
\bigskip
\heading II. Properties and hypotheses
\endheading
\bigskip
In this section, we give some properties and hypotheses for the
functions.

Now we consider the following basic Volterra integral equation
$$
\aligned x(t)=
  &x_0 (t)+\int_0^{t}f_1 (t, s, x(s))ds +\int_0^t \int_0^s f_2 (t,
    s, s_1 , x(s), x(s_1 ))ds_1 ds\\
  &+\sum_{i:\tau_i <t} G_1 (t, \tau_i , x(\tau_i ^- ))
   +\sum_{i:\tau_i <t}\sum_{j=1} ^{i-1} G_2 (t, \tau_i , \tau_j , x(\tau_i ^- ), x(\tau_j ^- ))\\
\endaligned
$$
$$
\aligned\sum
  &+\int_0^t \sum_{i:\sigma _i (s) <t}\sum_{j:\tau_j <t}
       g(t, s, \sigma _i (s), \tau _j , x(s), x(\sigma _i (s)^- ), x(\tau_j ^- ))ds\\
  &+\sum_{i:\sigma _i (t) <t}\sum_{j:\tau_j <t}
       G_3 (t, \sigma _i (t), \tau_j , x(\sigma _i (t)^- ), x(\tau _j ^-
       )).
\endaligned \leqno (2.1)
$$

Let $\rho_{ij}$, $i=1,2, \cdots , N_\sigma$, $j=1,2, \cdots ,
N_{\sigma ,i}$ be the solutions of equation
$$
t=\sigma_i (t) \leqno (2.2)
$$
Define $\sigma =\{\sigma _i (\cdot ):i=1,2, \cdots ,N_\sigma\}$,
and we assume that, for every $i=1,2,\cdots, N_\sigma$. equation
$(2.2)$ has a finite number of solutions. Let
$$
{\bold {\rho}}:=\{ \rho_{ij} :1\leq j\leq N_{\sigma, i} ,1\leq
i\leq N_\sigma \}. \leqno (2.3)
$$
We set
$$
{\bold \tau} :=\{\tau_0 , \tau_1 , \cdots ,\tau_{N_\tau}\} \leqno
(2.4)
$$
and
$$
{\bold I}:={\bold \tau} \cup {\bold \rho} \leqno (2.5)
$$

The space $C(0,T;R;{\bold I})$ is defined as the space of
real-valued functions $x(\cdot )$ which are bounded and continuous
on every open interval $(\alpha , \beta )$, with endpoints
$\alpha$ and $\beta$ in ${\bold I}$, and has limits $x(\alpha ^+
):=\lim_{t\longrightarrow \alpha^+ }x(t)$, $x(\alpha ^-
):=\lim_{t\longrightarrow \alpha^- }x(t)$ at every point $\alpha
\in {\bold I}\cap (0,T)$, as well as limits $x(0^+ ):=\lim_{t
\longrightarrow 0^+}x(t)$, $x(T^- ):=\lim_{t \longrightarrow T^-}
x(t)$.
\bigskip
Let $x(t)=\xi (t)$, $x(\tau_i ^- )=\eta_i$ and $x(\sigma_i (t)^-
)=\beta_i (t)$. Assume that the functions $f_1 , f_2 , G_1 , G_2 ,
G_3$ and $g$ satisfy the following conditions:

\noindent (H1) $f_1 (t,s,x(s))$ is continuous for $0\leq s \leq
t$, $x\in C(0,T;R;{\bold I})$, there exist $L_1
>0$ such that
$$
\vert f_1 (t,s,x_1 )-f_1 (t,s,x_2 )\vert \leq L_1 \vert x_1 -x_2
\vert .
$$

\noindent(H2) $f_2 (t, s, s_1 , x, y)$ is continuous for $0\leq
s_1 \leq s \leq t$, $x, y\in C(0,T;R;{\bold I})$, there exist
$L_{21}, L_{22}>0$ such that
$$
\vert f_2 (t,s, s_1 , x_1 ,y_1 )-f_2 (t,s,s_1 , x_2 ,y_2 )\vert
\leq L_{21} \vert x_1 -x_2 \vert +L_{22} \vert y_1 -y_2 \vert .
$$

\noindent (H3) $G_1 (t, \tau_i , x(\tau_i ^- ))$ is continuous
impulse function for $0<\tau_i <t<T$, $x \in R^n$, there exist
$L_{G_1} >0$ such that
$$
\vert G_1 (t, \tau_i , \eta_i^1 )-G_1 (t, \tau_i , \eta_i^2 )\vert
\leq L_{G_1} \vert \eta_i^1 -\eta_i^2 \vert .
$$

\noindent (H4) $G_2 (t, \tau_i , \tau_j , x(\tau_i ^- ), x(\tau_j
^- ))$ is continuous impulse function for $0<\tau_j <\tau_i <t<T$,
$x \in R^n$, there exist $L_{G_{21}}, L_{G_{22}} >0$ such that
$$
\vert G_2 (t, \tau_i , \tau_j , \eta_i^1 , \eta_j^1 )-G_2 (t,
\tau_i , \tau_j , \eta_i^2 , \eta_j^2 )\vert \leq L_{G_{21}} \vert
\eta_i^1 -\eta_i^2 \vert +L_{G_{22}} \vert \eta_j^1 -\eta_j^2
\vert .
$$

\noindent (H5) $G_3 (t, \sigma _i (t), \tau_j , x(\sigma _i (t)^-
), x(\tau _j ^- ))$ is continuous impulse function for $0<\tau_j
<\sigma_i (t) <t<T$, $x \in R^n$, there exist $L_{G_{31}},
L_{G_{32}}
>0$ such that
$$
\vert G_3 (t, \sigma _i (t), \tau_j , \beta_i^1 (t), \eta_j^1
)-G_3 (t, \sigma _i (t), \tau_j , \beta_i^2 (t), \eta_j^2 ) \vert
\leq L_{G_{31}} \vert \eta_j^1 -\eta_j^2 \vert +L_{G_{32}} \vert
\beta_i^1 (t)-\beta_i^2 (t)\vert .
$$

\noindent (H6) $g(t, s, \sigma _i (s), \tau _j , x(s), x(\sigma _i
(s)^- ), x(\tau_j ^- ))$ is continuous function for $0<\tau_j
<\sigma_i (s) <s<t<T$, $x\in C([0,T]:\tau )$, there exist $L_{g_1}
, L_{g_2} , L_{g_3} >0$ such that
$$
\aligned &\vert g(t, s, \sigma _i (s), \tau _j , \xi^1 (t),
\beta_i^1 (t), \eta_j^1 )-\vert g(t, s, \sigma _i (s), \tau _j ,
\xi^2 (t),
\beta_i^2 (t), \eta_j^2 )\vert \\
&\qquad \qquad \ \ \leq L_{g_1} \vert \xi^1 (t)-\xi^2 (t)\vert +
L_{g_2} \vert \eta_j^1 -\eta_j^2 \vert +L_{g_3} \vert \beta_i^1
(t) -\beta_i^2 (t)\vert .
\endaligned
$$
\noindent (H7) There is a positive number $h$ such that $\tau_i
-\tau_{i-1} \geq h $ for all $i=1,2,\cdots, N_\tau$;  $\sigma_j
(t) -\sigma_{j-1} (t)\geq h$ for all $j=1,2, \cdots , N_\sigma$,
and for all $t\in [0,T]$;  whenever $\sigma_j (s) <\tau_i$, for
some $j\in \{1,2, \cdots, N_\sigma \}$ and $i\in \{1,2, \cdots,
N_\tau \}$, then $\tau_i -\sigma _j (s)\geq h$.
\bigskip
\bigskip
\heading III. Existence and Uniqueness of solution of the state
equation
\endheading
\bigskip
In this section, we will show the existence and uniqueness of
solution for the nonlinear Volterra integral equation.

For each collection of impulse times, we seek a solution $(2.1)$
in the space $C(0,T;R;{\bold I})$ of real valued function $x(t)$
that are bounded and continuous on every open interval $(\alpha ,
\beta )$, with endpoints $\alpha$ and $\beta$ in ${\bold I}$ and
have limits $x(\alpha ^+ ):=\lim_{t\longrightarrow \alpha^+
}x(t)$, $x(\alpha ^- ):=\lim_{t\longrightarrow \alpha^- }x(t)$ at
every point $\alpha \in {\bold I}\cap (0,T)$, as well as limits
$x(0^+ ):=\lim_{t \longrightarrow 0^+}x(t)$, $x(T^- ):=\lim_{t
\longrightarrow T^-} x(t)$.

\bigskip
{\ \ \bf Theorem 3.1.} Suppose that hypotheses $(H1)\thicksim
(H7)$ are satisfied. Then for every $\tau$ and every $\sigma$,
$(2.1)$ has a unique solution $x(\cdot )$ in the space
$C(0,T;R;{\bold I} )$.
\bigskip
{\ \ \bf Proof.} We observe that the hybrid Volterra equation
$(2.1)$ implies the following impulsive conditions at the times
$\tau_k$ and $\rho_{kl}$:
$$
\aligned x(\tau_k ^+ )=
&x(\tau_k ^- )+G_1 (\tau_k , \tau_k ,
   x(\tau_k ^- ))+ \sum_{j=1} ^{k-1} G_2 (\tau_k , \tau_k , \tau_j ,
   x(\tau_k ^- ), x(\tau_j ^- ))\\
&+\int_0^{\tau_k} \sum_{i:\sigma_i (s)<\tau_k} g(\tau_k , s,
   \sigma_i (s),\tau_k ,x(s), x(\sigma_i (s)^- ), x(\tau_k ^- ))ds\\
&+\sum_{i:\sigma_i (\tau_k )<\tau_k} G_3 (\tau_k ,\sigma_i (\tau_k
   ),\tau_k ,x(\sigma_i (\tau_k )^- ), x(\tau_k ^- ))
\endaligned \leqno (3.1)
$$
$$
x(\rho_{kl} ^+ )=x(\rho_{kl} ^- )+\sum_{j:\tau_j <\rho_{kl}} G_3
(\rho_{kl} , \rho_{kl} ,\tau_j ,x(\rho_{kl} ^- ),x(\tau_j ^- ))
\leqno (3.2)
$$

The interval $[0,T]$ can be expressed as $[0,T]=\cup_{0\leq l\leq
M} [\alpha_l ,\alpha_{l+1}]$, for some positive integer $M$, with
each $\alpha_l \in {\bold I}$, and all the corresponding open
intervals satisfying $(\alpha_l , \alpha_{l+1} )\cap (\alpha_k ,
\alpha_{k+1} )= \emptyset$ for $k\neq l$, $(\alpha_l ,\alpha_{l+1}
)\cap {\bold I}=\emptyset$.

The solution of $(2.1)$ can be obtained inductively as follows:
equation $(2.1)$ can be solved on $[0,\alpha_1 )$ as a Volterra
integral equation(including double integral terms, for which
existence and uniqueness is provided by the argument in the
introduction), and therefore $x(\alpha_1 ^- )$ can  be determined.
Then $x(\alpha_1 ^+ )$ can be found by the impulsive conditions
$(3.1),(3.2)$ above. The equation $(2.1)$ becomes a Volterra
integral equation(with multiple integral terms, but no impulses)
for the restriction of $x(t)$ to the interval $(\alpha_1 ,\alpha_2
)$, and consequently $x(t)$ can be uniquely determined over the
interval $[0,\alpha_2 )$. Inductively, if $x(t)$ has been
determined over the interval $[0, \alpha_l )$, then $x(\alpha_l ^-
)$ can be determined, and $x(\alpha_l ^+)$ can be found from the
impulsive conditions $(3.1),(3.2)$, so that $(2.1)$ becomes a
Volterra integral equation(with multiple integral terms, but no
impulses) for the restriction of $x(t)$ to the interval $(\alpha_l
,\alpha_{l+1} )$, which is uniquely solvable on $(\alpha_l
,\alpha_{l+1} )$, therefore $x(t)$ will be known over the interval
$[0,\alpha_{l+1} )$. This completes the induction.
\bigskip
We consider the space $V=C(0,T;R;{\bold I})\times R^{N_\tau}
\times C(0,T;R;{\bold I})$ with a vector valued norm defined for
each $(\xi, \eta_i , \beta_i )$ in $V$ denoted by $\xi (t)=x(t)$,
$\eta_i =x(\tau_i^- )$, $\beta_i (t)=x(\sigma_i (t)^- )$ by
$$
\aligned
&\Vert \xi \Vert_\mu =\sup_{0\leq t \leq T}e^{-\mu t}
\vert \xi (t)\vert\\
&\Vert \eta \Vert_\mu =\max_{1\leq i\leq N} e^{-\mu \tau_i}
\vert \eta_i \vert , \ \ \ \ \ \ \ \ \ \ \ \ \ \eta =[\eta_i :0\leq i \leq N_\tau ]\\
&\Vert \beta \Vert_\mu =\max_{1\leq i\leq N}\sup_{0\leq t \leq T}
e^{-\mu \sigma_i (t)} \vert \beta_i (t)\vert , \ \ \ \ \ \ \beta
=[\beta_j :0\leq j\leq N_\sigma ]
\endaligned \leqno (3.3)
$$

We define an operator $S$ on $V$ by
$$
S(\xi , \eta ,\beta )= \pmatrix
 S_c (\xi, \eta, \beta )\\
 [S_d (\xi, \eta, \beta )]_l\\
 [S_m (\xi, \eta, \beta )]_p ,
\endpmatrix \leqno (3.4)
$$
where $S_c :V \longrightarrow C(0,T;R;{\bold I})$, $S_d
:V\longrightarrow R^{N_\tau}$ and $S_m : V\longrightarrow
C(0,T;R;{\bold I})$;
$$
\aligned &S_c (\xi, \eta, \beta )(t)\\
  &=x_0 (t)+\int_0^{t}f_1 (t, s, \xi(s))ds +\int_0^t \int_0^s f_2 (t,
    s, s_1 , \xi(s), \xi(s_1 ))ds_1 ds\\
  &\ \ \ \ +\sum_{i:\tau_i <t} G_1 (t, \tau_i , \eta_i )
   +\sum_{i:\tau_i <t}\sum_{j=1} ^{i-1} G_2 (t, \tau_i , \tau_j , \eta_i , \eta_j )\\
\endaligned
$$
$$
\aligned
  &\ \ \ \ +\int_0^t \sum_{i:\sigma _i (s) <t}\sum_{j:\tau_j <t}
       g(t, s, \sigma _i (s), \tau _j , \xi (s), \beta _i (s), \eta_j )ds\\
  &\ \ \ \ +\sum_{i:\sigma _i (t) <t}\sum_{j:\tau_j <t}
       G_3 (t, \sigma _i (t), \tau_j , \beta _i (t), \eta_j ),
\endaligned \leqno (3.5)
$$
$$
\aligned &[S_d (\xi, \eta, \beta )]_l \\
  &=x_0 (\tau_l )+\int_0^{\tau_l}f_1 (\tau_l , s, \xi(s))ds +\int_0^{\tau_l} \int_0^s f_2 (\tau_l ,
    s, s_1 , \xi(s), \xi(s_1 ))ds_1 ds\\
  &\ \ \ \ +\sum_{i<l} G_1 (\tau_l , \tau_i , \eta_i )
   +\sum_{i<l}\sum_{j=1} ^{i-1} G_2 (\tau_l , \tau_i , \tau_j , \eta_i , \eta_j )\\
  &\ \ \ \ +\int_0^{\tau_l} \sum_{i:\sigma _i (s) <\tau_l}\sum_{j=1}^{l-1}
       g(\tau_l , s, \sigma _i (s), \tau _j , \xi (s), \beta _i (s), \eta_j )ds\\
  &\ \ \ \ +\sum_{i:\sigma _i (\tau_l ) <\tau_l}\sum_{j=1}^{l-1}
       G_3 (\tau_l , \sigma _i (\tau_l ), \tau_j , \beta _i (\tau_l ), \eta_j ),
\endaligned \leqno (3.6)
$$
$$
\aligned &[S_m (\xi, \eta, \beta )]_p (t)\\
  &=x_0 (\sigma_p (t))+\int_0^{\sigma_p (t)}f_1 (\sigma_p (t), s, \xi(s))ds \\
  &\ \ \ \ +\int_0^{\sigma_p (t)} \int_0^s f_2 (\sigma_p (t),
    s, s_1 , \xi(s), \xi(s_1 ))ds_1 ds\\
  &\ \ \ \ +\sum_{i:\tau_i <\sigma_p (t)} G_1 (\sigma_p (t), \tau_i , \eta_i )
   +\sum_{i:\tau_i < \sigma_p (t)}\sum_{j=1} ^{i-1} G_2 (\sigma_p (t), \tau_i , \tau_j , \eta_i , \eta_j )\\
  &\ \ \ \ +\int_0^{\sigma_p (t)} \sum_{i:\sigma _i (s) <\sigma_p (t)}\sum_{j:\tau_j<\sigma_p (t)}
       g(\sigma_p (t), s, \sigma _i (s), \tau _j , \xi (s), \beta _i (s), \eta_j )ds\\
  &\ \ \ \ +\sum_{i:\sigma _i (\sigma_p (t)) <\sigma_p (t)}\sum_{j:\tau_j<\sigma_p (t)}
       G_3 (\sigma_p (t), \sigma _i (\sigma_p (t)), \tau_j , \beta _i (\sigma_p (t)), \eta_j ).
\endaligned \leqno (3.7)
$$

We shall call that $S_c$ is continuous component of $S$, $S_d$ is
discrete component of $S$ and $S_m$ is mixed component of $S$.
\bigskip
{\ \ \bf Lemma 3.2.} The solution of equation $(2.1)$ is
equivalent to the problem of finding a fixed point of the operator
$S$ defined in $(3.5)\sim (3.7)$ above.
\bigskip
{\ \ \bf Proof.} It is clear that, if $x(t)$ is a solution of
$(2.1)$ in the space $C(0,T;R; {\bold I})$, then it is plain that
the triple $(\xi^* , \eta^* , \beta^* )$, defined by
$$
\xi^* (t):=x(t), \  \ \eta_i ^* :=x(\tau_i ^- ), \ \ \beta_j ^*
(t):=x(\sigma_j (t)^- )
$$
is a fixed point of the operator $S$ in the space $V$.

Conversely, suppose that $(\xi^* , \eta^* , \beta^* )\in V$ is a
fixed point of $S$. Let the intervals $[\alpha_l ,\alpha_{l+1} ]$
be as in the proof of Theorem 3.1. It follows from $(3.5)$ that
$\xi^* (t)$ is a solution of $(2.1)$ over the time-interval
$[0,\alpha_1 )$. If $\alpha_1 \in {\bold \tau}$, it follows from
$(3.6)$ that $\eta_1 ^* =\xi^* (\alpha_1 ^- )$, and then it
follows from $(3.5)$ that $\xi^* (\alpha_1 ^+ )$ satisfies the
impulsive condition $(3.1)$ at $t=\alpha_1$. If $\alpha_1 \in
{\bold \rho}$, it follows from $(3.7)$ that $\beta_1 ^* (\alpha_1
)=\xi^* (\sigma_1 (\alpha_1 )^- )$, and then it follows from
$(3.5)$ that $\xi^* (\sigma_1 (\alpha_1 )^+ )$ satisfies the
impulsive condition $(3.2)$ at $t=\alpha_1$. Inductively, if
$\xi^* (t)$ solves $(2.1)$ over $[0,\alpha_l )$ and $\xi^*
(\alpha_l ^+ )$ satisfies the appropriate impulsive condition,
$(3.1)$ or $(3.2)$, at $t=\alpha_l$, then we shall show that
$\xi^* (t)$ also solves $(2.1)$ over $[0,\alpha_{l+1} )$ and
$\xi^* (\alpha_{l+1} ^+ )$ satisfies the appropriate impulsive
condition, $(3.1)$ or $(3.2)$, at $t=\alpha_{l+1}$. Suppose that
the points $\alpha_0 , \alpha_1 , \cdots , \alpha_l$ correspond to
$\{ \tau_\lambda : 0\leq \lambda \leq \lambda_0 \}\cup \{
\sigma_\mu : 1\leq \sigma_\mu \leq \mu_0 \}$. Then it follows from
$(3.5)$ that $\xi^* (t)$ solves $(2.1)$ over $[0,\alpha_{l+1} )$,
since the equation obtained from $(2.1)$ over the interval
$[\alpha_l ,\alpha_{l+1} )$ as in the proof of Theorem 3.1,
contains no impulses in $(\alpha_l , \alpha_{l+1} )$. If
$\alpha_{l+1} \in {\bold \tau}$, then it follows from $(3.6)$ that
$\xi^* (\alpha_{l+1} ^- )=\eta_{\lambda_0 +1} ^*$, and then it
follows from $(3.5)$ that, in case $\alpha_{l+1} <T, \ \xi^*
(\alpha_{l+1} ^+ )$ satisfies the impulsive condition $(3.1)$ at
$t=\alpha_{l+1} =\tau_{\lambda_0 +1}$. If $\alpha_{l+1} \in {\bold
\rho}$, then it follows from $(3.6)$ that $\xi^* (\alpha_{l+1} ^-
)=\beta_j ^* (\sigma_j (\alpha_{l+1} )^- )$ for some $j\geq
\mu_0$, and then it follows from $(3.5)$ that, in case
$\alpha_{l+1} <T, \ \xi^* (\alpha_{l+1} ^+ )$ satisfies the
impulsive condition $(3.2)$ at $t=\rho_{jk} =\alpha_{l+1}$ for
some $k$. The inclusion is complete.
\bigskip
{\ \ \bf Lemma 3.3.} If $\mu>0$ is sufficiently large,
then for any $(\xi^1 ,\eta^1 ,\beta^1 ),(\xi^2 ,\eta^2 ,\beta^2 )$
of $V$, we have
$$\aligned
&{\text (i)}\ \Vert S_c (\xi^1 ,\eta^1 ,\beta^1 )-S_c (\xi^2
,\eta^2
,\beta^2 )\Vert _\mu \\
&\qquad \qquad \qquad \qquad \qquad \leq a_{11} \Vert \xi^1 -\xi^2
\Vert _\mu +a_{12}\Vert \eta^1 -\eta^2 \Vert_ \mu +a_{13}\Vert
\beta^1 -\beta^2 \Vert _\mu,
\endaligned
$$
$$ \aligned
&{\text (ii)}\ \Vert S_d (\xi^1 ,\eta^1 ,\beta^1 )-S_d (\xi^2
,\eta^2
,\beta^2 )\Vert _\mu \\
&\qquad \qquad \qquad \qquad \qquad \leq a_{21} \Vert \xi^1 -\xi^2
\Vert _\mu +a_{22}\Vert \eta^1 -\eta^2 \Vert_ \mu +a_{23}\Vert
\beta^1 -\beta^2 \Vert _\mu,
\endaligned
$$
$$ \aligned
&{\text (iii)}\ \Vert S_m (\xi^1 ,\eta^1 ,\beta^1 )-S_m (\xi^2
,\eta^2
,\beta^2 )\Vert _\mu \\
&\qquad \qquad \qquad \qquad \qquad \leq a_{31} \Vert \xi^1 -\xi^2
\Vert _\mu +a_{32}\Vert \eta^1 -\eta^2 \Vert_ \mu +a_{33}\Vert
\beta^1 -\beta^2 \Vert _\mu,
\endaligned
$$
where all constants $a_{ij}$ are nonnegative elements of the
matrix $A:=(a_{ij})_{1\leq i,j \leq 3}$
\bigskip
{\ \ \bf Proof.} We will show that the operator $S$ satisfies the
contractive condition. It is enough to have $\lim_{\mu
\longrightarrow \infty}D=0$, $\lim_{\mu \longrightarrow
\infty}T=0$ and $\lim_{\mu \longrightarrow \infty}S=0$ for
$D=det(A)$, $T=tr(A)$ and
$S=a_{11}a_{22}+a_{22}a_{33}+a_{33}a_{11}-a_{13}a_{31}-a_{23}a_{32}-a_{12}a_{21}$
(see  Appendix $I$).

Let $L_2 =\max \{ L_{21}, L_{22} \}$, $L_{G_2}=\max \{L_{G_{21}},
L_{G_{22}} \}$, $L_g =\max \{L_{g_1}, L_{g_2},L_{g_3} \}$ and
$L_{G_3} =\max \{L_{G_{31}}, L_{G_{32}}\}$. Then we have
$$
\aligned &{\text (i)}\  e^{-\mu t} \vert S_c (\xi^1 ,\eta^1
,\beta^1
   )(t)-S_c (\xi^2 ,\eta^2 ,\beta^2 )(t) \vert\\
&\qquad \ \ \leq e^{-\mu t} \biggl{\{} \int_0^t \vert f_1 (t,
    s, \xi ^1(s))-f_1 (t, s, \xi ^2 (s))\vert ds \\
&\qquad \ \ \ \ \ +\int_0^t \int_0^s \vert f_2 (t, s, s_1 , \xi ^1
(s), \xi ^1
    (s_1 ))-f_2 (t, s, s_1 , \xi ^2 (s), \xi ^2 (s_1 ))\vert ds_1 ds\\
&\qquad \ \ \ \ \ +\sum_{i:\tau_i <t} \vert G_1 (t, \tau_i ,
\eta_i ^1 )-G_1
     (t, \tau_i , \eta_i ^2 )\vert \\
&\qquad \ \ \ \ \ +\sum_{i:\tau_i <t}\sum_{j=1} ^{i-1} \vert G_2
(t, \tau_i ,
       \tau_j , \eta_i ^1 , \eta_j ^1 )
            -G_2 (t, \tau_i , \tau_j , \eta_i ^2 , \eta_j ^2 )\vert\\
&\qquad \ \ \ \ \ +\int_0^t \sum_{i:\sigma _i (s) <t}
\sum_{j:\tau_j <t}
          \vert g(t, s, \sigma _i (s), \tau _j , \xi^1 (s), \beta _i ^1 (s), \eta_j ^1
          )\\
&\qquad \qquad \qquad \ \ -g(t, s, \sigma _i (s), \tau _j , \xi^2 (s), \beta _i ^2 (s), \eta_j ^2 )\vert ds\\
&\qquad \ \ \ \ \ +\sum_{i:\sigma _i (t) <t}\sum_{j:\tau_j <t}
          \vert G_3 (t, \sigma _i (t), \tau_j , \beta _i ^1 (t), \eta_j ^1 )
            -G_3 (t, \sigma _i (t), \tau_j , \beta _i ^2 (t), \eta_j ^2
            )\vert \biggr{\}}.
\endaligned
$$
Hence, by the assumption $(H1)\sim(H7)$,
$$
\aligned
&\Vert S_c (\xi^1 ,\eta^1 ,\beta^1 )-S_c (\xi^2 ,\eta^2
    ,\beta^2 ) \Vert _\mu\\
&\leq \biggl{\{} L_1 \frac{1-e^{-\mu t}}{\mu}+L_2 (\frac{t}{\mu}
   -\frac{te^{-\mu t}}{\mu})+L_g (N_t )^2 \frac{1-e^{-\mu
    t}}{\mu}\biggr{\}}\Vert \xi^1 -\xi^2 \Vert _\mu\\
&\ \ +\biggl{\{} L_{G_1} (1+\frac{1-e^{-\mu (N_t -1)h}}{e^{\mu h}
           -1})+L_{G_2} (N_t -1+(N_t -1)(N_t -2)+N_t (N_t -1))e^{-\mu h}\\
&\qquad \ \ +L_g N_t (1+\frac{1-e^{-\mu (N_t -1)h}}{e^{\mu h}-1}T+
      L_{G_3} N_t (1+\frac{1-e^{-\mu (N_t -1)h}}{e^{\mu h}-1} )\biggr{\}}\Vert \eta^1 -\eta^2 \Vert _\mu\\
&\ \ +\biggl{\{}L_{G_3} N_t (\frac{1-e^{-\mu N_t h}}{1-e^{-\mu
h}})+L_g TN_t (\frac{1-e^{-\mu N_t h}}{1-e^{-\mu h}})\biggr{\}}
\Vert \beta^1 -\beta^2 \Vert _\mu .
\endaligned
$$
$$
\aligned &{\text (ii)}\  e^{-\mu \tau_l} \vert S_d (\xi^1 ,\eta^1
   ,\beta^1
   )(\tau_l )-S_d (\xi^2 ,\eta^2 ,\beta^2 )(\tau_l ) \vert\\
&\qquad \ \ \leq e^{-\mu \tau_l} \biggl{\{} \int_0^{\tau_l} \vert
f_1 (\tau_l ,s, \xi ^1(s))-f_1 (\tau_l , s, \xi ^2 (s))\vert ds \\
&\qquad \ \ \ \ \ +\int_0^{\tau_l} \int_0^s \vert f_2 (\tau_l , s,
s_1 , \xi ^1 (s), \xi ^1
    (s_1 ))-f_2 (\tau_l , s, s_1 , \xi ^2 (s), \xi ^2 (s_1 ))\vert ds_1 ds\\
&\qquad  \ \ \ \ \ +\sum_{i<l} \vert G_1 (\tau_l , \tau_i , \eta_i
^1 )-G_1 (\tau_l , \tau_i , \eta_i ^2 )\vert \\
\endaligned
$$
$$
\aligned
&\qquad \ \ \ \ \ +\sum_{i<l}\sum_{j=1} ^{i-1} \vert G_2
(\tau_l , \tau_i ,\tau_j , \eta_i ^1 , \eta_j ^1 )
            -G_2 (\tau_l , \tau_i , \tau_j , \eta_i ^2 , \eta_j ^2 )\vert\\
&\qquad \ \ \ \ \ +\int_0^{\tau_l} \sum_{i:\sigma _i (s) <\tau_l}
+\sum_{j=1} ^{l-1} \vert g(\tau_l , s, \sigma _i (s), \tau _j ,
\xi^1 (s), \beta _i ^1 (s), \eta_j ^1 )\\
&\quad \qquad \qquad \ \ -g(\tau_l , s, \sigma _i (s), \tau _j , \xi^2 (s), \beta _i ^2 (s), \eta_j ^2 )\vert ds\\
&\qquad \ \ \ \ \ +\sum_{i:\sigma _i (\tau_l )
<\tau_l}\sum_{j=1}^{l-1}
          \vert G_3 (\tau_l , \sigma _i (\tau_l ), \tau_j , \beta _i ^1 (\tau_l ), \eta_j ^1 )
            -G_3 (\tau_l , \sigma _i (\tau_l ), \tau_j , \beta _i ^2 (\tau_l ), \eta_j
            ^2
            )\vert \biggr{\}}.
\endaligned
$$
\bigskip
Hence, by the assumption $(H1)\sim(H7)$,
$$
\aligned &\Vert S_d (\xi^1 ,\eta^1 ,\beta^1 )-S_d (\xi^2 ,\eta^2
    ,\beta^2 ) \Vert _\mu\\
&\leq \biggl{\{} L_1 \frac{1-e^{-\mu \tau_l}}{\mu}+L_2
(\frac{\tau_l}{\mu}
   -\frac{\tau_l e^{-\mu \tau_l}}{\mu})+L_g (l-1)N_{\tau_l} \frac{1-e^{-\mu
    \tau_l}}{\mu}\biggr{\}}\Vert \xi^1 -\xi^2 \Vert _\mu\\
&\ \ +\biggl{\{} L_{G_1} \frac{1-e^{-\mu (l -1)h}}{e^{\mu h}
           -1}+L_{G_2} \frac{l(l-1)e^{-\mu h}}{2}+L_{G_2}\frac{1-e^{-\mu (l-2)h}}{(e^{\mu h}
           -1)^2}\\
&\qquad \qquad \qquad +L_g N_{\tau_l} \tau_l \frac{1-e^{-\mu (l
-1)h}}{e^{\mu h}-1}+
      L_{G_3} N_{\tau_l} \frac{1-e^{-\mu (l -1)h}}{e^{\mu h}-1}\biggr{\}}\Vert \eta^1 -\eta^2 \Vert _\mu\\
&\ \ +\biggl{\{}L_{g} (l-1)\tau_l e^{-\mu
h}\frac{1-e^{-\mu(N_{\tau_l} -1)h}}{1-e^{-\mu h}}+L_{G_3}
(l-1)e^{-\mu h} \frac{1-e^{-\mu (N_{\tau_l}-1)h}}{1-e^{-\mu
h}})\biggr{\}} \Vert \beta^1 -\beta^2 \Vert _\mu .
\endaligned
$$
$$
\aligned &{\text (iii)}\  e^{-\mu t} \vert [S_m (\xi^1 ,\eta^1
    ,\beta^1
   )]_p (t)-[S_m (\xi^2 ,\eta^2 ,\beta^2 )]_p (t) \vert\\
&\qquad \ \ \leq e^{-\mu t} \biggl{\{} \int_0^{\sigma_p (t)} \vert
f_1 (\sigma_p (t),s, \xi ^1(s))-f_1 (\sigma_p (t), s, \xi ^2 (s))\vert ds \\
&\qquad  \ \ \ \ \ +\int_0^{\sigma_p (t)} \int_0^s \vert f_2
(\sigma_p (t), s, s_1 , \xi ^1 (s), \xi ^1
    (s_1 ))-f_2 (\sigma_p (t), s, s_1 , \xi ^2 (s), \xi ^2 (s_1 ))\vert ds_1 ds\\
&\qquad  \ \ \ \ \ +\sum_{i:\tau_i <\sigma _p (t)} \vert G_1
(\sigma_p (t), \tau_i , \eta_i
^1 )-G_1 (\sigma_p (t), \tau_i , \eta_i ^2 )\vert \\
&\qquad  \ \ \ \ \ +\sum_{i:\tau_i <\sigma_p (t)}\sum_{j=1} ^{i-1}
\vert G_2 (\sigma_p (t), \tau_i ,\tau_j , \eta_i ^1 , \eta_j ^1 )
            -G_2 (\sigma_p (t), \tau_i , \tau_j , \eta_i ^2 , \eta_j ^2 )\vert\\
&\qquad \ \ \ \ \ +\int_0^{\sigma_p (t)} \sum_{i:\sigma _i (s)
<\sigma_p (t)} \sum_{j:\tau_j <\sigma_p (t)} \vert g(\sigma_p (t),
s, \sigma _i (s), \tau _j ,
\xi^1 (s), \beta _i ^1 (s), \eta_j ^1 )\\
&\quad \qquad \qquad \qquad \qquad -g(\sigma_p (t), s, \sigma _i (s), \tau _j , \xi^2 (s), \beta _i ^2 (s), \eta_j ^2 )\vert ds\\
\endaligned
$$
$$
\aligned
&\qquad \ \ \ \ \ +\sum_{i:\sigma _i (\sigma_p (t))
<\sigma_p (t)}\sum_{j:\tau_j <\sigma_p (t)}
          \vert G_3 (\sigma_p (t), \sigma _i (\sigma_p (t)), \tau_j , \beta _i ^1 (\sigma_p (t)), \eta_j ^1
          )\\
&\quad \qquad \qquad \qquad \qquad -G_3 (\sigma_p (t), \sigma _i
(\sigma_p (t)), \tau_j , \beta _i ^2 (\sigma_p (t)), \eta_j ^2
)\vert \biggr{\}}.
\endaligned
$$
Assume that $\sigma_{i+1} (t)-\sigma_i (s)\geq h$ whenever $s\leq
\sigma_i (t)$. Then, by the assumptions $(H1)\sim(H7)$,
$$
\aligned &\Vert S_m (\xi^1 ,\eta^1 ,\beta^1 )-S_m (\xi^2 ,\eta^2
    ,\beta^2 ) \Vert _\mu\\
&\leq \biggl{\{} L_1 \frac{1-e^{-\mu \sigma_p (t)}}{\mu}+L_2
(\frac{\sigma_p (t)}{\mu}
   -\frac{\sigma_p (t)e^{-\mu \sigma_p (t)}}{\mu})+L_g (N_{\sigma_p (t)})^2 \frac{1-e^{-\mu
    \sigma_p (t)}}{\mu}\biggr{\}}\Vert \xi^1 -\xi^2 \Vert _\mu\\
&\ \ +\biggl{\{} L_{G_1} \frac{1-e^{-\mu hN_{\sigma_p
(t)}}}{1-e^{-\mu h}}+L_{G_2} (N_{\sigma_p (t)}+\frac{N_{\sigma_p
(t)} (N_{\sigma_p (t)-1)e^{-\mu h}}}{2}
   -\frac{1-e^{-\mu hN_{\sigma_p (t)}}}{1-e^{-\mu h}})\\
&\qquad \qquad \qquad +L_{G_2}(N_{\sigma_p (t)}e^{-\mu h}
+\frac{N_{\sigma_p (t)}(N_{\sigma_p (t)} -1)e^{-\mu
h}}{2}-\frac{1-e^{-\mu
N_{\sigma_p (t)} h}}{e^{\mu h} -1})\\
&\qquad \qquad \qquad +L_g \sigma_p (t)N_{\sigma_p
(t)}\frac{1-e^{-\mu hN_{\sigma_p (t)}}}{1-e^{\mu h}}+
      L_{G_3} N_{\sigma_p (t)} \frac{1-e^{-\mu hN_{\sigma_p (t)}}}{1-e^{-\mu h}}\biggr{\}}\Vert \eta^1 -\eta^2 \Vert _\mu\\
&\ \ +\biggl{\{}L_{g} \sigma_p (t)N_{\sigma_p
(t)}\frac{1-e^{-\mu(N_{\sigma_p (t)} -1)h}}{e^{\mu h}-1}+L_{G_3}
N_{\sigma_p (t)}\frac{1-e^{-\mu (N_{\sigma_p (t)}-1)h}}{e^{\mu
h}-1})\biggr{\}} \Vert \beta^1 -\beta^2 \Vert _\mu .
\endaligned
$$

If we take $a_{ij}, \ 1\leq i,j\leq 3$, in the above result
inequalities as same as $(i)$, $(ii)$ and $(iii)$ of this lemma,
then $(a_{11}, a_{21}, a_{22}, a_{23}, a_{31}, a_{33})
\longrightarrow (0,0,0,0,0,0)$ as $\mu \longrightarrow \infty $,
and $a_{12}, a_{13}, a_{32}$ remains bounded as $\mu
\longrightarrow \infty$ ; consequently $(D,T,S)\longrightarrow 0$
as $\mu \longrightarrow\infty$. That is, for $\mu$ sufficiently
large, the operator $S$ will be satisfied the contractive
condition.
\bigskip
We consider the following iterative scheme for the solution of
$(2.1)$: $(\xi(0), \eta(0), \beta(0))$ is an arbitrary element of
$V$; for $k=0, 1, 2, \cdots$, $(\xi(k+1), \eta(k+1), \beta(k+1))$
is defined by
$$
\aligned &\xi_{(k+1)} (t)\\
  &=x_0 (t)+\int_0^{t}f_1 (t, s, \xi_{(k)}(s))ds +\int_0^t \int_0^s f_2 (t,
    s, s_1 , \xi_{(k)}(s), \xi_{(k)} (s_1 ))ds_1 ds\\
  &\quad +\sum_{i:\tau_i <t} G_1 (t, \tau_i , \eta_{(k),i} )
   +\sum_{i:\tau_i <t}\sum_{j=1} ^{i-1} G_2 (t, \tau_i , \tau_j , \eta_{(k),i} , \eta_{(k),j} )\\
  &\quad+\int_0^t \sum_{i:\sigma _i (s) <t}\sum_{j:\tau_j <t}
       g(t, s, \sigma _i (s), \tau _j , \xi_{(k)} (s), \beta_{(k),i} (s), \eta_{(k),j} )ds\\
  &\quad+\sum_{i:\sigma _i (t) <t}\sum_{j:\tau_j <t}
       G_3 (t, \sigma _i (t), \tau_j , \beta _{(k),i} (t), \eta_{(k),j} ),
\endaligned \leqno (3.8)
$$
$$
\aligned &\eta_{(k+1),l}\\
  &=x_0 (\tau_l )+\int_0^{\tau_l}f_1 (\tau_l , s, \xi_{(k)}(s))ds +\int_0^{\tau_l} \int_0^s f_2 (\tau_l ,
    s, s_1 , \xi_{(k)}(s), \xi_{(k)}(s_1 ))ds_1 ds\\
  &\quad+\sum_{i<l} G_1 (\tau_l , \tau_i , \eta_{(k),i} )
   +\sum_{i<l}\sum_{j=1} ^{i-1} G_2 (\tau_l , \tau_i , \tau_j , \eta_{(k),i} , \eta_{(k),j} )\\
  &\quad+\int_0^{\tau_l} \sum_{i:\sigma _i (s) <\tau_l}\sum_{j=1}^{l-1}
       g(\tau_l , s, \sigma _i (s), \tau _j , \xi_{(k)} (s), \beta_{(k),i} (s), \eta_{(k),j} )ds\\
  &\quad+\sum_{i:\sigma _i (\tau_l ) <\tau_l}\sum_{j=1}^{l-1}
       G_3 (\tau_l , \sigma _i (\tau_l ), \tau_j , \beta_{(k),i} (\tau_l ), \eta_{(k),j} ),
\endaligned \leqno (3.9)
$$
$$
\aligned &\beta_{(k+1),p} (t)\\
  &=x_0 (\sigma_p (t))+\int_0^{\sigma_p (t)}f_1 (\sigma_p (t), s, \xi_{(k)}(s))ds \\
  &\quad+\int_0^{\sigma_p (t)} \int_0^s f_2 (\sigma_p (t),
    s, s_1 , \xi_{(k)}(s), \xi_{(k)}(s_1 ))ds_1 ds\\
  &\quad+\sum_{i:\tau_i <\sigma_p (t)} G_1 (\sigma_p (t), \tau_i , \eta_{(k),i} )
   +\sum_{i:\tau_i < \sigma_p (t)}\sum_{j=1} ^{i-1} G_2 (\sigma_p (t), \tau_i , \tau_j , \eta_{(k),i} , \eta_{(k),j} )\\
  &\quad+\int_0^{\sigma_p (t)} \sum_{i:\sigma _i (s) <\sigma_p (t)}\sum_{j:\tau_j<\sigma_p (t)}
       g(\sigma_p (t), s, \sigma _i (s), \tau _j , \xi_{(k)}(s), \beta _{(k),i} (s), \eta_{(k),j} )ds\\
  &\quad+\sum_{i:\sigma _i (\sigma_p (t)) <\sigma_p (t)}\sum_{j:\tau_j<\sigma_p (t)}
       G_3 (\sigma_p (t), \sigma _i (\sigma_p (t)), \tau_j , \beta _{(k),i} (\sigma_p (t)), \eta_{(k),j} ).
\endaligned \leqno (3.10)
$$
Then we have:
\bigskip
{\ \ \bf Theorem 3.4.} Suppose that $(H1)\sim(H7)$ are satisfied.
Then the iterative method defined by $(3.8)\sim (3.10)$ above
converges to $(x(\cdot ), \ (x(\tau_1^- ),x(\tau_2^- ),x(\tau_2^-
),$ $\cdots, x(\tau_{N_\tau} ^- ))$, $x(\sigma_1 (\cdot )^-
),x(\sigma_2 (\cdot )^- ),$ $ \cdots, x(\sigma_{N_\sigma} (\cdot
)^- ))$, as $k \longrightarrow \infty$; the convergence of
$\xi_{(k)}$ to $x(\cdot )$ is uniform on each closed interval
$[\alpha , \beta ]$ that does not contain points of $\bold I$ in
the interior, in the sense that, if we define the restrictions
$\xi_{(k)} ^{\alpha \beta} =\xi_{(k)}, \ x^{\alpha \beta} (\cdot
)=x $ to the interval $[\alpha, \beta]$ by $\xi_{(k)} ^{\alpha
\beta} (t)=\xi_{(k)}(t)$
 for $t\in [\alpha ,\beta )$, $\xi_{(k)}^{\alpha \beta}(
 \beta )=\xi_{(k)}(\beta ^- )$, $x^{\alpha \beta}(t)=x(t)$ for
 $t\in [\alpha , \beta )$, $x^{\alpha \beta} (\beta )=x(\beta ^-
 )$, then $\xi_{(k)} ^{\alpha \beta}\longrightarrow x^{\alpha
 \beta}$ uniformly on $[\alpha ,\beta ]$.
\bigskip
{\ \ \bf Proof.} By Lemma 3.3, if $\mu$ is sufficiently large, the
operator $S$ is a contraction with respect to the vector-valued
norm $\Vert \cdot \Vert_\mu$ on $V$. The contraction property with
respect to the vector valued norm $\Vert \cdot \Vert_\mu$ means
that, for all $x,y \in V$, we have $\Vert Sx -Sy\Vert_\mu \leq A
\Vert x-y\Vert_\mu$, where the $3\times 3$ real matrix $A$ has
$\lim_{\mu \longrightarrow \infty}D=0$, $\lim_{\mu \longrightarrow
\infty}T=0$, $\lim_{\mu \longrightarrow \infty}S=0$. Consequently,
the iterates of $S$, with arbitrary initial data, converges to the
unique fixed point of $S$ in the topology induced on $V$ by the
vector valued norm $\Vert \cdot \Vert_\mu$; this is a well known
extension of the standard Banach fixed point theorem to the case
of a vector valued metric, and the proof proceeds as in the
standard case. Convergence with respect to uniform convergence on
each $[\tau_{i-1} , \tau_i ]$. The fixed point of S gives the
solution of $(2.1)$ by Lemma 3.2.
\bigskip
{\ \ \bf Remark.} It follows from $(3.5)\sim (3.7)$ that
$\eta_{(k+1),l}=\xi_{(k+1)} (\tau_l^- ), \
\beta_{(k+1),l}=\xi_{(k)} (\sigma_l (t)^- )$ for all $k\geq 0$, so
that, for $k\geq 1$, $(3.5)\sim (3.7)$ can also be written in the
form
$$
\aligned \xi_{(k+1)} (t)=
  &x_0 (t)+\int_0^{t}f_1 (t, s, \xi_{(k)}(s))ds +\int_0^t \int_0^s f_2 (t,
    s, s_1 , \xi_{(k)}(s), \xi_{(k)} (s_1 ))ds_1 ds\\
  &+\sum_{i:\tau_i <t} G_1 (t, \tau_i , \xi_{(k)}(\tau^- ))
   +\sum_{i:\tau_i <t}\sum_{j=1} ^{i-1} G_2 (t, \tau_i , \tau_j , \xi_{(k)}(\tau_i ^- ) , \xi_{(k)}(\tau^- ))\\
  &+\int_0^t \sum_{i:\sigma _i (s) <t}\sum_{j:\tau_j <t}
       g(t, s, \sigma _i (s), \tau _j , \xi_{(k)} (s), \xi_{(k)} (s), \xi_{(k)}(\tau_j^- ))ds\\
  &+\sum_{i:\sigma _i (t) <t}\sum_{j:\tau_j <t}
       G_3 (t, \sigma _i (t), \tau_j , \xi _{(k)}(\sigma_i (t)^- ), \xi_{(k)}(\tau_j ^- )
       ).
\endaligned
$$
Of  course, if $\eta_{(0)}$ is chosen as
$\eta_{(0),i}=\xi_{(0)}(\tau_i ^- )$, $\beta_{(0)}$ is chosen as
$\beta_{(0),i}=\xi_{(0)}(\sigma_i (t)^- )$, then $(2.1)$ holds for
all $k=0,1,2, \cdots \cdots$
\bigskip
{\ \ \bf Appendix I. Condition for a $3\times 3$ matrix to be
contractive.} We shall derive the necessary and sufficient
conditions for a $3\times 3$ matrix $A=[a_{ij}]_{1\leq i, j\leq
3}$ to be contractive, in the sense that all eigenvalues of $A$
will have modulus strictly less than 1.

The characteristic polynomial of $A$ is
$$
f(\lambda )=\lambda^3 -T\lambda^2 +S\lambda -D \leqno (A.1)
$$
where $T=tr(A)=a_{11} +a_{22} +a_{33}$, the trace of $A$, $D=\det
(A)$, the determinant of $A$, and
$$
\aligned
S&={\biggl \vert}{\aligned & a_{11} \ \  a_{12}\\
                          & a_{21} \ \  a_{22}\endaligned}{\biggr \vert}
  +{\biggl \vert}{\aligned & a_{11} \ \  a_{13}\\
                 & a_{31} \ \  a_{33}\endaligned}{\biggr \vert}
 +{\biggl \vert}{\aligned & a_{22} \ \  a_{23}\\
                 & a_{32} \ \  a_{33}\endaligned}{\biggr \vert}\\
  &=a_{11} a_{22} +a_{22} a_{33}
  +a_{33} a_{11}-a_{13} a_{31}-a_{23} a_{32}-a_{12} a_{21} .
\endaligned \leqno (A.2)
$$

We recall, from the theory of the Routh-Hurwitz stability
criteria, that a cubic polynomial $\varphi (\omega )=p_0 \omega^3
+p_1 \omega^2 +p_2 \omega +p_3$, with $p_0 >0$, will have all its
roots in the open left-half plane (on the plane of complex
numbers) if and only if the following conditions are satisfied:
$$
p_1 >0;\ \ \ \ \  p_3 >0;\ \ \ \ \  p_1 p_2 >p_0 p_3 . \leqno
(A.3)
$$

The transformation $\lambda =\frac{\omega +1}{\omega -1}$ maps
$Re(\omega )<0$ onto $\vert \lambda \vert <1$. By using this
transformation into $(A.1)$, we find that $A$ will have all its
eigenvalues in the interior of the unit circle on the complex
plane if and only if the polynomial
$$
\aligned g(\omega )=
     &[1-T+S-D]\omega^3 +[3-T-S+3D]\omega^2\\
     &+[3+T-S-3D]\omega +[1-T+S+D]
\endaligned \leqno (A.4)
$$
\noindent has all its roots in the half-plane $Re(\omega )<0$.
After some straightforward algebra, we find
$$
p_1 p_2 -p_0 p_3 =8-4S-8D^2 +4DT. \leqno (A.5)
$$

Consequently, assuming $p_0\equiv 1-T+S-D>0$, the matrix $A$ will
be contractive if and only if
$$
\aligned &3-S+3D-T>0 ;\\
         &1+T+S+D>0;\\
         &8-4S-8D^2 +4DT >0
\endaligned  \leqno (A. 6)
$$
Further, we note that a sufficient condition for $A$ to be
contractive is that the quantities $S, T, D$ should be
sufficiently small.

In case the elements $a_{ij}$ of the matrix $A$ depend on a
parameter $\mu$ (as will be the case in the application to the
mappings defined in Section 3 of this paper, a sufficient
condition for $A$ to be contractive is that
$(S,D,T)\longrightarrow (0,0,0)$ as $\mu \longrightarrow \infty .$
\bigskip
\bigskip
\centerline{\bf  References}
\bigskip
\item {$ \text{[ 1 ]}$} D. Bainov \& P.S. Simeonov(1989). Systems with impulse effect,
                         Ellis Horwood Ltd, Chichester.
\item {$ \text{[ 2 ]}$} S.A. Belbas(1999). Iterative schemes for optimal control of Volterra integral equations,
                         {\it Nonlinear Analysis} 37, 57-79.
\item {$ \text{[ 3 ]}$} S.A. Belbas \& W.E. Schmidt(2004).  Optimal control of Volterra equations with impulses,
         {\it Applied Math.} \& {\it Computation}.
                                         {\it Fuzzy Sets and System},  24, 301-317.
\item {$ \text{[ 4 ]}$} C. Corduneanu(1991). Integral equations and applications, Cambridge University Press,
                       Cambridge.
\item {$ \text{[ 5 ]}$} D. Guo(1992). Impulsive integral equations in Banach spaces and applications,
                     {\it J. Appl. Math. $\&$ Stochastic Anal. } 5, 111-122.
\item {$ \text{[ 6 ]}$} D. Guo(1993). Nonlinear impulsive Volterra integral equations in Banach spaces and applications,
                     {\it J. Appl. Math. $\&$ Stochastic Anal. }
                     6, 35-48
\item {$ \text{[ 7 ]}$} S. Hong(2004). Solvability of nonlinear
           impulsive Volterra integral inclusions and functional differential inclusions, {\it J.
           Math. Anal. $\&$ Appl.} 295, 331-340.
\item {$ \text{[ 8 ]}$} J. Kondo(1991). Integral equations, Kodansa, Tokyo, and Clarendon Press, Oxford.
\item {$ \text{[ 9 ]}$} J.S. Park, Y.C. Kwun, C.W. Han \& S.Y. Kim(2004).
       Existence and Uniqueness of fuzzy solutions for the nonlinear fuzzy
      integro-differential equation on $E_N^n$.
       {\it International J. of Fuzzy Logic and Intelligent Systems},  4(1), 40-44.
\item {$ \text{[ 10 ]}$} J.S. Park \& Y.C. Kwun(2004).
           Existence results for semilinear integro-differential equations
           using integral contractors.
           {\it Far East J. Math. Sci.}
           12(1), 65-72.
\item {$ \text{[ 11 ]}$} B. Yan(1997). On $L_{loc}^p$ -solutions of
       nonlinear impulsive Volterra integral equations in Banach spaces,
       {\it SUT J. Math.} 33, 121-137.
\item {}
\item {} {\it E-mail address} : *sbelbas \@bama.ua.edu
\item {} {\it E-mail address} : **parkjs\@cue.ac.kr
\enddocument